\documentclass{scrartcl}

\usepackage[T1]{fontenc}
\usepackage[latin1]{inputenc}

\pdfoutput=1

\usepackage{amssymb}
\usepackage{amsfonts}
\usepackage{amsmath}
\usepackage{wasysym}

\usepackage{epsfig}
\usepackage{epic}
\usepackage{eepic}
\usepackage{graphics}
\usepackage{graphicx}
\usepackage{color}
\usepackage{longtable}

\usepackage{amsthm}

\newcommand{\ds}{\displaystyle}
\newcommand{\ts}{\textstyle}

\newcommand{\var}{\operatorname{var}}
\newcommand{\spa}{\operatorname{span}}
\newcommand{\supp}{\operatorname{supp}}

\renewcommand{\Re}{\operatorname{Re}}

\newcommand{\Nn}{{\mathbb N}}

\newcommand{\Rr}{{\mathbb R}}

\newcommand{\vph}{\varphi}
\newcommand{\eps}{\varepsilon}
\newcommand{\vvar}{\mathrm{var}}

\newcommand{\Jj}{\mathbf{J}}
\newcommand{\ccc}{\mathbf{c}}

\newcommand{\pia}{\Pi_n}
\newcommand{\pib}{\Pi_{n}^m}

\def\be{\begin{equation}}
\def\ee{\end{equation}}
\def\ba{\begin{align*}}
\def\ea{\end{align*}}

\theoremstyle{plain}
\newtheorem{Thm}{Theorem}[section]

\newtheorem{Lem}[Thm]{Lemma}

\theoremstyle{remark}
\newtheorem{Rem}[Thm]{Remark}
\newtheorem{Exa}[Thm]{Example}

\theoremstyle{definition}
\newtheorem{Def}[Thm]{Definition}

\usepackage{tikz}
\usetikzlibrary{calc,fadings,decorations.pathreplacing}
\usetikzlibrary{matrix,arrows}

\begin{document}

\title{An orthogonal polynomial analogue of the Landau-Pollak-Slepian time-frequency analysis}
\author{Wolfgang Erb
\thanks{Institute of Mathematics, University of Lübeck,
Ratzeburger Allee 160, 23562 Lübeck, Germany. erb@math.uni-luebeck.de}}

\date{15.03.2012}

\maketitle

\begin{abstract}
The aim of this article is to present a time-frequency theory for orthogonal polynomials on the interval $[-1,1]$ 
that runs parallel to the time-frequency analysis of bandlimited functions developed by Landau, Pollak and Slepian. For this purpose, 
the spectral decomposition of a particular compact time-frequency-operator is studied. This decomposition and its eigenvalues are
closely related to the theory of orthogonal polynomials. 
Results from both theories, the theory of orthogonal polynomials and the Landau-Pollak-Slepian theory, can be used to prove localization and
approximation properties of the corresponding eigenfunctions. Finally, an uncertainty principle is proven that reflects the limitation of coupled time and frequency locatability. 
\end{abstract}


\section{Introduction}

In the beginning of the 1960s, Landau, Pollak and Slepian developed a remarkable theory on the time-frequency analysis of band-limited functions.
In a series of papers (\cite{LandauPollak1961}, \cite{LandauPollak1962}, \cite{LandauWidom1980}, 
\cite{Slepian1964}, \cite{Slepian1978}, \cite{SlepianPollak1961}) they studied the interplay between the two projection operators $P_A$ and $P_B$ defined on the Hilbert
space $L^2(\Rr)$ for two intervals $A,B \subset \Rr$ by
\[ P_A f := \chi_A f, \qquad \widehat{P_B f} := \chi_B \hat{f}, \qquad f \in L^2(\Rr).\]
They analyzed the composition $P_B P_A P_B$ and its spectrum and found that the eigenfunctions of the compact self-adjoint operator $P_B P_A P_B$ are well-known special functions: the prolate spheroidal wave functions. Using these particular eigenfunctions as a basis for the band-limited functions in $L^2(\Rr)$ on the other hand, they were able to prove a series of interesting results concerning the approximate concentration of functions in the time and the frequency domain, as well as an uncertainty principle involving a lower bound for the angle between the vectors $P_A f$ and $P_B f$. An overview of these results can be found in the articles \cite{Landau1985}, \cite{Slepian1983} and the book \cite[Section 2.9]{DymMcKean}.

Later on, the Landau-Pollak-Slepian-theory was extended to a variety of different settings. Among others, there exist analogies on the unit circle \cite{Slepian1983}, on discrete groups
\cite{Grünbaum1981} and on symmetric spaces like the unit sphere \cite{GrünbaumLonghiPerlstadt1982}, \cite{SimonsDahlenWieczorek2006}. Various generalizations of this theory can be formulated, for instance by considering eigenfunctions of particular Sturm-Liouville differential equations \cite{WangZhang2010} or using reproducing kernel Hilbert spaces \cite{Zayed2007}. Particularly interesting for this article is the fact that there exists also an extension of this theory to orthogonal polynomials defined on subsets of the real line \cite{Perlstadt1986}. 

The aim of this paper is to present a time-frequency analysis for orthogonal polynomials on the interval $[-1,1]$ 
that runs parallel to the Landau-Pollak-Slepian theory described in \cite{Perlstadt1986}. For the frequency localization
of a function $f$ in the weighted $L^2$-space $L^2([-1,1],w)$ we will use, as in \cite{Perlstadt1986}, an operator $P_n^m$ that 
projects the function $f$ onto a finite dimensional 
polynomial space $\pib$. However, in contrast to the theory outlined above, we will not use a projection operator $P_A$ to describe 
the space localization of $f$. Instead, we will consider the multiplication operator $M_x$ defined by multiplying the function $f$ with the variable $x$. 

Compared to the projection operator $P_A$, the usage of the multiplication operator $M_x$ leads to a time-frequency analysis in which the localization of $f$ at the boundary points $x=1$ and $x=-1$ of the interval $[-1,1]$ plays an important role. For a normalized function $f \in L^2([-1,1],w)$, the mean value $\eps(f) = \langle M_x f, f\rangle_w$ is located in the interval $(-1,1)$. The closer $\eps(f)$ gets to $1$ or $-1$, the more the $L^2$-mass of $f$ is concentrated at $x = 1$ or $x = -1$, respectively. Therefore, the mean value $\eps(f)$ can be considered as a measure on how well the function $f$ is localized at the boundary points $x=1$ or $x=-1$. Particularly this property of $\eps(f)$ implies the possibility to construct polynomials in $\pib$ that are optimally localized at the boundary of $[-1,1]$
(see \cite{Erb2012}, \cite{Freud1986}, \cite{Rauhut2005}). 

The principal examination object for the time-frequency analysis in this paper is the finite dimensional self-adjoint operator $P_n^m M_x P_n^m$ in combination with its eigenvalues $x_{n,k}^m$, $1 \leq k \leq n-m+1$, and corresponding eigenfunctions $\psi_{n,k}^m$. One of the main advantages of the operator $M_x$ in place of $P_A$ is the fact that the spectral decomposition of $P_n^m M_x P_n^m$ is closely linked to the theory of orthogonal polynomials. This relation makes it possible to use a very large repertoire of techniques and results from the theory of orthogonal polynomials to analyse the properties of the spectral decomposition of $P_n^m M_x P_n^m$. In the spectral Theorem \ref{Thm-spectraldecomposition}, we will see that the eigenvalues of $P_n^m M_x P_n^m$ are precisely the roots of the associated orthogonal 
polynomials $p_{n-m+1}(x,m)$. Also the eigenfunctions can be stated explicitly. In the case $m = 0$, they correspond to the fundamental polynomials of Lagrange interpolation. 

A second advantage of using the operator $M_x$ consists in the fact that the value $\eps(f)$ represents also the expectation value of the $L^2$-density $f$. The density $f$ can be considered as localized at the expected value $\eps(f)$ if the variance $\var(f)$ is small. Therefore, we can investigate the localization properties of the eigenfunctions $\psi_{n,k}^m$ of $P_n^m M_x P_n^m$ by considering the variances $\var(\psi_{n,k}^m)$. In order to show that the functionals $\var(\psi_{n,k}^m)$ are small when $n$ is large, we will use results of Nevai, Zhang and Totik \cite{Nevai}, \cite{NevaiZhangTotik1991} on uniform subexponential growth. 
The major result in this context is Theorem \ref{Theorem-localizedeigenfunctions}. It states
that if the weight function $w$ of the space $L^2([-1,1],w)$ is in a particular subclass of the Nevai class $M(0,1)$, 
then the variance of the eigenfunctions $\psi_{n,k}^m$ tend to zero as $n \to \infty$. 

In Section \ref{section-approximationfunctions}, we will analyse how the decomposition of a bandlimited function $f \in \pib$ in the single eigenfunctions $\psi_{n,k}^m$ can be used to approximate functions that are localized at a point or a subinterval of $[-1,1]$. In this case, not all the eigenfunctions $\psi_{n,k}^m$ are needed to approximate the function $f$, but just those that are situated in the region in which $f$ is concentrated. In Theorem \ref{Theorem-localizedapproximationepsilonconcentrated} and \ref{Theorem-localizedapproximationvariance} we will give simple error estimates for such approximations if the function $f$ is localized in a certain area or at a particular point of the interval $[-1,1]$, respectively. 

Finally, we will prove an uncertainty relation for orthogonal polynomials involving the operators $M_x$ and $P_n^m$. This relation can be considered as an extension of the angular uncertainty principle in the Landau-Pollak-Slepian theory. For a normalized function $f \subset L^2([-1,1],w)$, the determining quantities of the uncertainty relation are the norm $\|P_n^m f\|_w$ and again the mean value $\eps(f)$. 
The norm $\|P_n^m f\|_w$ gives a measure on how well the function $f$ is concentrated in the polynomial subspace $\Pi_n^m$. On the other hand the value $\eps(f)$ can be seen as a measure of the localization of $f$ at the boundary points $x=-1$ and $x = 1$. The main result in the last section is Theorem \ref{Theorem-uncertainty} claiming that for a normalized function $f \in L^2([-1,1],w)$ it is impossible that $\|P_n^m f\|_w$ and $|\eps(f)|$ are both close to $1$. In particular, this result implies that if $|\eps(f)|$ is too close to $1$, $f$ cannot be a polynomial in $\Pi_n^m$.

\section{The spectral decomposition of $P_n^m M_x P_n^m$}

We consider the Hilbert space $L^2([-1,1],w)$ with the inner product
\[ \langle f,g \rangle_{w} := \int_{-1}^1 f(x) \overline{g(x)} w(x) dx,\]
and a positive weight function $w$ having finite moments $\int_{-1}^1 x^n w(x) dx$, $n \in \Nn$. By $\{p_l\}_{l=0}^\infty$, we denote the family of polynomials $p_l$ of degree $l$ that are orthonormal on $[-1,1]$ with respect to the inner product
$\langle \cdot, \cdot \rangle_w$. Further, we assume that the polynomials $p_l$ are normalized such that the coefficient of the monomial $x^l$ is positive. Then, the family $\{p_l\}_{l=0}^\infty$ defines a complete orthonormal set in the Hilbert space $L^2([-1,1],w)$ (cf. \cite[Section 2.2]{Szegö}). By $\pia$, we denote the polynomial space spanned by the polynomials $p_l$ up to degree $n$, and by $\pib$ the polynomial wavelet space spanned by the polynomials $p_l$, $m \leq l \leq n$.

For a normalized function $f \in L^2([-1,1],w)$, $\|f\|_w = 1$, we define the mean value $\eps(f)$ and the variance $\vvar(f)$ by
\begin{align}
\eps(f) &:= \int_{-1}^1 x |f(x)|^2 w(x) dx, \label{equation-meanvalue} \\
\vvar(f) & := \int_{-1}^1 (x - \eps(f))^2 |f(x)|^2 w(x) dx = \int_{-1}^1 x^2 |f(x)|^2 w(x) dx - \eps(f)^2. \label{equation-variance}
\end{align}

We are now going to introduce a time-frequency analysis for functions $f \in L^2([-1,1],w)$ based on the following two operators:
\begin{align}
 (M_x f)(x)  &:= x f(x), \\
 (P_n^m f)(x) &:= \sum_{l=m}^n \langle f, p_l \rangle_w\, p_l(x).
\end{align}
If $m=0$, we write $P_n$ instead of $P_n^0$. The multiplication operator $M_x$ as well as the orthogonal projection $P_n^m$ onto $\pib$ are both self-adjoint and bounded operators
on the Hilbert space $L^2([-1,1],w)$. Therefore,
also the composition
\begin{equation}
P_n^m M_x P_n^m
\end{equation}
is a bounded and self-adjoint operator on $L^2([-1,1],w)$. Moreover, since $P_n^m$ is compact, $P_n^m M_x P_n^m$ is also a compact
operator. Hence, by the Hilbert-Schmidt theorem the spectrum of the operator $P_n^m M_x P_n^m$ is discrete (it is even finite) 
and the eigenfunctions form an orthogonal basis of $L^2([-1,1],w)$ (cf. \cite[Theorem VI.16]{ReedSimon1}). 
The subsequent Theorem \ref{Thm-spectraldecomposition} will illustrate that the eigenvalues and eigenfunctions of $P_n^m M_x P_n^m$ are well-known in the literature.

For a description of the spectral decompostion of $P_n^m M_x P_n^m$, we need first of all the notion of associated polynomials.
We know that the orthonormal polynomials $p_l$ satisfy the three-term recurrence relation
(cf. \cite[Section 1.3.2]{Gautschi})
\begin{align} \label{equation-recursionorthonormal}
b_{l+1} p_{l+1}(x) &= (x - a_l) p_l(x) - b_{l} p_{l-1}(x), \quad l=0,1,2,3, \ldots \\
 p_{-1}(x) &= 0, \qquad p_0(x) = \frac{1}{b_0}, \notag
\end{align}
with coefficients $a_l \in \Rr$ and $b_l > 0$. For $m \in \Nn$, the associated polynomials $p_l(x,m)$ on the interval $[-1,1]$ are
then defined by the shifted recurrence relation (see \cite[Section 1.3.4]{Gautschi}, \cite[Section 2.10]{Ismail})
\begin{align} \label{equation-recursionassociatedsymmetric}
b_{m+l+1}\, p_{l+1}(x,m) &= (x - a_{m+l})\, p_l(x,m) - b_{m+l}\, p_{l-1}(x,m),
\quad l=0,1,2, \ldots , \\
 p_{-1}(x,m) &= 0, \qquad p_0(x,m) = 1. \notag
\end{align}
For $m =0$, we have the identity, $p_l(x,0) = b_0\, p_l(x)$.
The polynomials $p_l(x)$ and $p_l(x,m)$
can be described with help of the symmetric Jacobi matrix $\Jj_n^m$, $0 \leq m \leq n $, defined by
\begin{equation} \label{equation-Jacobimatrix}
\Jj_n^m = \left(\begin{array}{cccccc}
a_m & b_{m+1} &  0 &  0 & \cdots & 0 \\
b_{m+1} &  a_{m+1} &  b_{m+2} &  0 & \cdots & 0 \\
0 &  b_{m+2}  & a_{m+2}  & b_{m+3} &  \ddots & \vdots \\
\vdots & \ddots & \ddots & \ddots & \ddots & 0\\
0 &  \cdots &  0&  b_{n-2} & a_{n-1} & b_{n-1} \\
0 &  \cdots &  \cdots & 0 & b_{n-1}& a_n
\end{array}\right).\end{equation}
If $m=0$, we write $\Jj_n$ instead of $\Jj_n^0$. Then, in view of the three-term recurrence formulas
\eqref{equation-recursionassociatedsymmetric}, the polynomials $p_l$ and $p_l(x,m)$ can be written as (cf. \cite[Theorem 2.2.4]{Ismail})
\begin{align}
p_l(x) &= \frac{1}{b_0} \det (x \mathbf{1}_l - \Jj_{l-1}), \label{equation-relation3termJacobimatrix} \\
p_l(x,m) &= \det(x \mathbf{1}_{l} - \Jj_{m+l-1}^m ), \label{equation-relation3termJacobimatrixassociated}
\end{align}
where $\mathbf{1}_{l}$ denotes the $l$-dimensional identity matrix. We can now explicitly state the spectral decomposition
of the operator $P_n^m M_x P_n^m$.

\begin{Thm} \label{Thm-spectraldecomposition}
The operator $P_n^m M_x P_n^m$ on $L^2([-1,1],w)$ has the spectral decomposition
\begin{equation}
P_n^m M_x P_n^m f = \sum_{k=1}^{n-m+1} x_{n,k}^m \langle f, \psi_{n,k}^m \rangle_w \psi_{n,k}^m.
\end{equation}
For $m \geq 1$, the eigenvalues $x_{n,k}^m$ denote the $n-m+1$ roots of the associated polynomial $p_{n-m+1}(x,m)$ and the eigenfunctions
$\psi_{n,k}^m$ have the explicit form
\begin{equation} \label{equation-slepianfunctionsexplicit}
\psi_{n,k}^m(x) = \kappa_{n,k}^m  \frac{ b_{n+1}p_{n+1}(x)
p_{n-m}(x_{n,k}^m,m)+ b_m p_{m-1}(x)}{x - x_{n,k}^{m}},
\end{equation}
with the normalizing constant
\begin{equation} \label{equation-normalizingconstant1}
\kappa_{n,k}^m := \Big( \sum_{l=m}^{n} p_{l-m}(x_{n,k}^m,m)^2 \Big)^{-\frac{1}{2}}.
\end{equation}
For $m = 0$, the eigenvalues $x_{n,k}$ correspond to the $n+1$ roots of the polynomial
$p_{n+1}(x)$ and the eigenfunctions $\psi_{n,k}$ correspond, up to a normalizing factor, to the fundamental polynomials of Lagrange interpolation, i.e.
\begin{equation} \label{equation-slepianfunctionsexplicit0}
\psi_{n,k}(x) = \kappa_{n,k}  p_n(x_{n,k}) b_{n+1} \frac{ p_{n+1}(x)}{x - x_{n,k}},
\end{equation}
where
\begin{equation}
\kappa_{n,k} := \Big( \sum_{l=0}^{n} p_{l}(x_{n,k})^2 \Big)^{-\frac{1}{2}}.
\end{equation}
\end{Thm}

\begin{proof}
We consider the projection $P_n^m f$ of the function $f$ onto the subspace $\pib$ in terms of the expansion $P_n^m f = \sum_{l=m}^n c_l p_l$ with the 
coefficients $c_l = \langle f, p_l \rangle_w$. Using the three term recurrence relation \eqref{equation-recursionassociatedsymmetric} it is straightforward
to show (see \cite[Lemma 2.7]{Erb2012}) that the mean value $\eps(P_n^m f)$ of $P_n^m f$ can be written as
\begin{equation} \label{equation-characterizationofeps} \langle P_n^m M_x P_n^m f,f \rangle_w = \langle M_x P_n^m f, P_n^m f \rangle_w =  \eps(P_n^m f) = \ccc^H \Jj_{n}^m \ccc,  \end{equation}
where $\ccc^H$ denotes the conjugate transpose of the vector $\ccc = (c_m, \ldots, c_n)^T$. 
Thus, the eigenvalues of $P_n^m M_x P_n^m$ in $\pib \subset L^2([-1,1],w)$ correspond to the eigenvalues of the Jacobi matrix $\Jj_n^m$.
On the other hand, by equation \eqref{equation-relation3termJacobimatrixassociated} the eigenvalues of $\Jj_n^m$ are exactly the
roots $x_{n,k}^m$, $k = 1, \ldots n-m+1$, of the associated polynomial $p_{n-m+1}(x,m)$. The eigenvector $\ccc_k$ corresponding to
the root $x_{n,k}^m$ is simple and can be computed via the three-term recursion formula \eqref{equation-recursionassociatedsymmetric} as
\begin{equation}  \label{equation-slepianfunctionscoefficients}
\ccc_k = \big(1, p_1(x_{n,k}^m,m), \ldots, p_{n-m}(x_{n,k}^m,m)\big)^T.
\end{equation}
The corresponding normalized eigenfunction $\psi_{n,k}^m$ of $P_n^m M_x P_n^m$ can then be written as
\begin{equation} \label{equation-slepianfunctionsfourierseries}
\psi_{n,k}^m(x) = \kappa_{n,k}^m \sum_{l=m}^{n} p_{l-m}(x_{n,k}^m,m) p_{l}(x),
\end{equation}
with the normalizing constant $\kappa_{n,k}^m$ given in \eqref{equation-normalizingconstant1}. By an alteration of the classical Christoffel-Darboux formula (see \cite[Lemma 3.1]{Erb2012}), the eigenfunctions $\psi_{n,k}^m$ for $m \geq 1$ have the explicit form
\[\psi_{n,k}^m(x) = \kappa_{n,k}^m  \frac{ b_{n+1}p_{n+1}(x)
p_{n-m}(x_{n,k}^m,m)+ b_m p_{m-1}(x)}{x - x_{n,k}^{m}}.\]
For $m=0$, we get directly by the Christoffel-Darboux formula (see \cite[Chapter 1, Theorem 4.5]{Chihara}) that 
\begin{equation*}
\psi_{n,k}(x) = \kappa_{n,k} b_{n+1}   \frac{ p_n(x_{n,k}) p_{n+1}(x)}{x - x_{n,k}}.
\end{equation*}
\end{proof}

\begin{Rem}
In the literature, the spectral Theorem \ref{Thm-spectraldecomposition} is well-known for the case $m = 0$ (cf. \cite[Lemma 8.4]{BreuerLastSimon2010}
and \cite[Proposition 1.3.1]{Simon2011}). For the more general case $m \geq 0$, an equivalent representation of Theorem \ref{Thm-spectraldecomposition} is the
eigenvalue decomposition $\Jj_n^m \ccc_k = x_{n,k}^m \ccc_k$ of the matrix $\Jj_n^m$ (see \cite[Section 1.3]{Gautschi}). To the best of the authors knowledge, the explicit formulas
\eqref{equation-slepianfunctionsexplicit} of the eigenfunctions $\psi_{n,k}^m$, $m \geq 1$, can be considered as novel. 
\end{Rem}


\begin{Rem}
The eigenfunctions $\{\psi_{n,k}^m\}_{k=1}^{n-m+1}$ of the operator $P_n^m M_x P_n^m$ form an orthonormal basis of the polynomial space $\pib$. Hence, we can expand 
polynomials $P \in \pib$ as
\[P(x) = \sum_{k=1}^{n-m+1} \langle P, \psi_{n,k}^m \rangle_w \psi_{n,k}^m(x).\]
In the case $m = 0$ the functions $\psi_{n,k}$ correspond to the fundamental polynomials of Lagrange interpolation and can be described through the 
Christoffel-Darboux kernel (see \cite[(1.1.9)]{Mhaskar} and formula \eqref{equation-slepianfunctionsfourierseries}). The functions $\psi_{n,k}$ are used in 
\cite{FischerPrestin1997} and \cite{FischerThemistoclakis2002} as particular orthogonal scaling functions in a wavelet decomposition of a function $f \in L^2([-1,1],w)$. 
If $m \geq 1$, the construction of the wavelet basis functions in these two papers differs however from the eigenfunctions $\psi_{n,k}^m$ considered in this article.  
For a general overview on polynomial frames and polynomial wavelet decompositions. we further refer to the articles \cite{FilbirMhaskarPrestin2009}, \cite{MhaskarPrestin2005} and the book \cite{Mhaskar}.
\end{Rem}

\begin{Rem} It was specified in the introduction that the mean value $\eps(f)$ can be interpreted as a measure on how localized the function $f$ is on the boundary points
$x = 1$ and $x=-1$ of $[-1,1]$. In the following, we will say that a function $f$ is localized at $x=1$ or $x=-1$ if the mean value $\eps(f)$ approaches $1$ or $-1$, respectively. 
For a polynomial $P \in \pib$, the mean value $\eps(P)$ can be written as $\eps(P) = \langle P_n^m M_x P_n^m P, P \rangle_w$. 
Precisely this mean value $\eps(P)$ was used in \cite{ErbDiss} and \cite{Erb2012} to construct polynomials in $\pia$ and $\pib$ that 
are optimally space localized at the boundary points $x = 1$ and $x = -1 $ of the interval $[-1,1]$. These optimal polynomials are exactly
the eigenfunctions $\psi_{n,\max}^m$ and $\psi_{n,\min}^m$ in Theorem \ref{Thm-spectraldecomposition} corresponding to the largest and the smallest eigenvalue of the operator $P_n^m M_x P_n^m$. 
By \eqref{equation-characterizationofeps}, we have for the largest eigenvalue of $P_n^m M_x P_n^m$ the relation
\[x_{n,\max}^m = \max_{P \in \pib, \|P\|_w = 1} \langle P_n^m M_x P_n^m P, P \rangle_w = \max_{\ccc^H \ccc = 1} \ccc^H \Jj_n^m \ccc.\]
This characterization is thoroughly used in \cite{Freud1986} to get estimates for the largest zero of orthogonal polynomials.

Taking a step further, we can also consider the orthogonal complement $\pib \ominus \spa \{\psi_{n,\max}^m\}$ of $\psi_{n,\max}^m$ in $\pib$. 
Then, the spectral Theorem \ref{Thm-spectraldecomposition} says that the polynomial in $\pib \ominus \spa \{\psi_{n,\max}^m\}$
that is best localized at $x = 1$ is the eigenfunction $\psi_{n,\max-1}^m$ corresponding to the second largest eigenvalue $x_{n,\max-1}^m$ of 
$P_n^m M_x P_n^m$. Hence, repeating this argumentation, Theorem \ref{Thm-spectraldecomposition} produces a chain of elementary orthonomal basis functions $\psi_{n,k}^m$
in which the $k$-th. element is worse concentrated at $x = 1$ than the $(k+1)$-th. element $\psi_{n,k+1}^m$ but better than the $(k-1)$-th. element $\psi_{n,k-1}^m$. The measure
of the corresponding localization is given by the mean value $\eps(\psi_{n,k}^m) = x_{n,k}^m$. 
\end{Rem}

\begin{figure}
  \begin{minipage}{0.5\textwidth}
  \centering
  $\psi_{24,25}(x)$. \\
  \includegraphics[width=\textwidth]{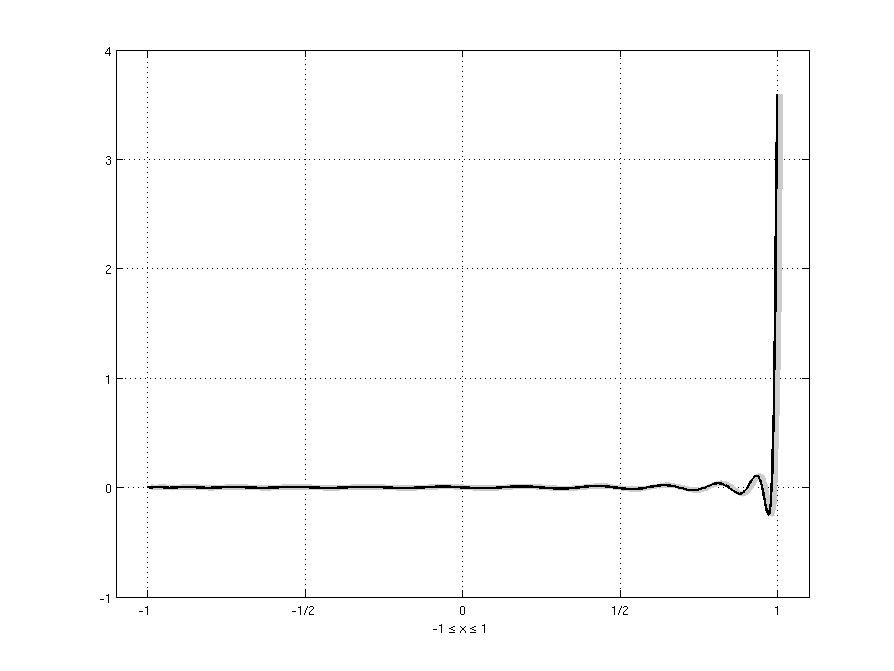}\\
  $\psi_{24,15}(x)$. \\
  \includegraphics[width=\textwidth]{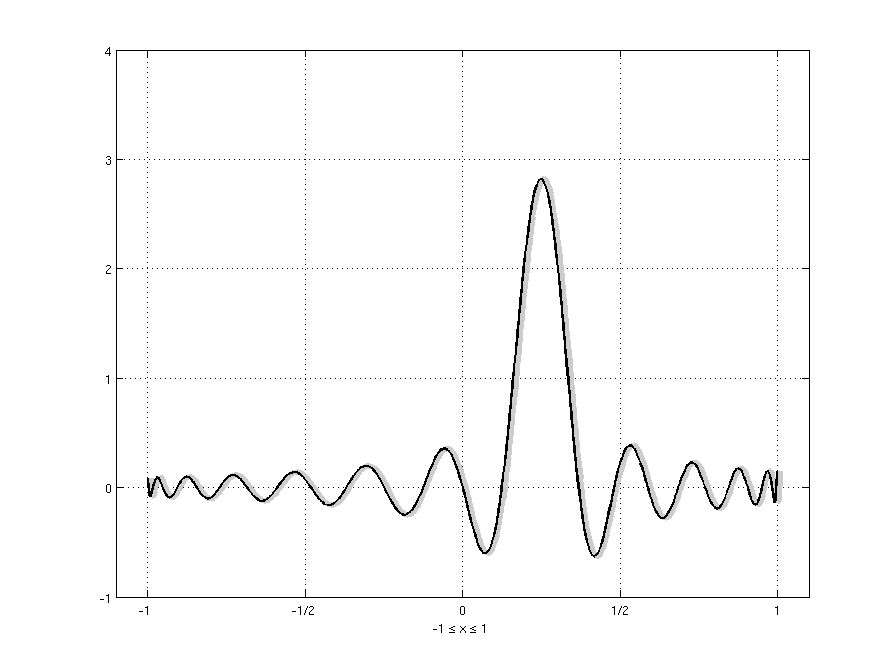}\\
  $\psi_{24,5}(x)$. \\
  \includegraphics[width=\textwidth]{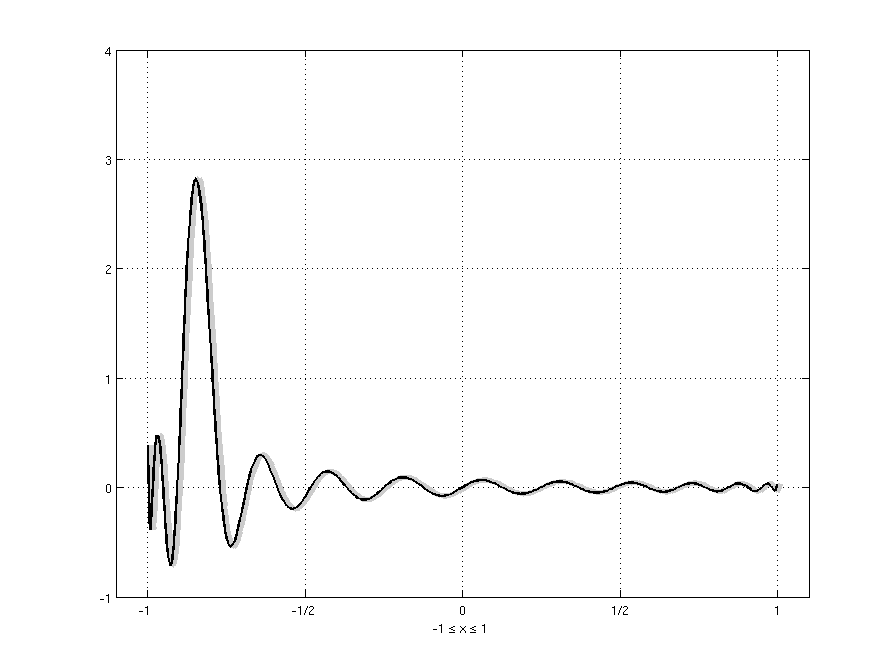}
  \end{minipage}\hfill
  \begin{minipage}{0.5\textwidth}
  \centering
  $\psi_{32,25}^8(x)$. \\
  \includegraphics[width=\textwidth]{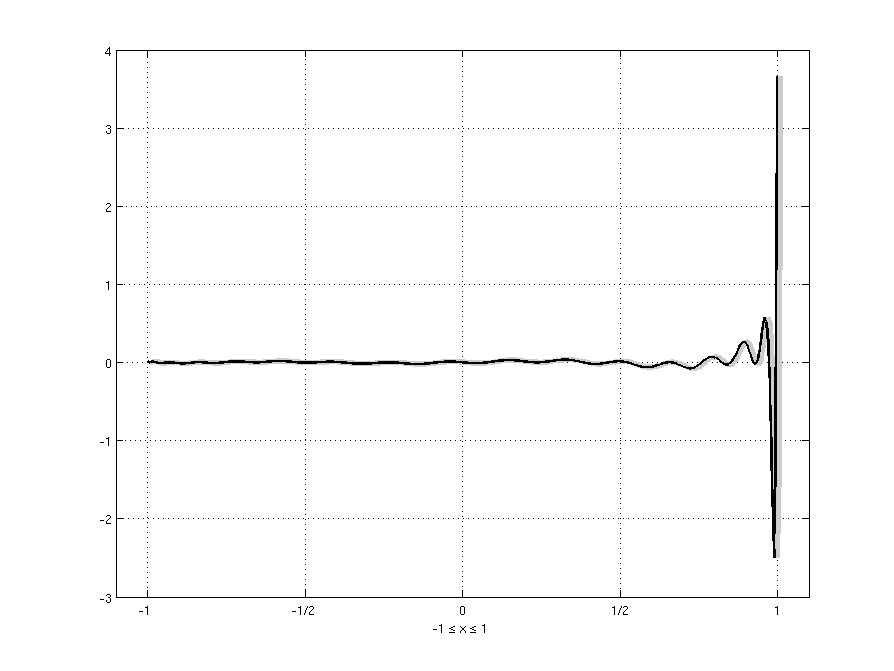}\\
  $\psi_{32,15}^8(x)$. \\
  \includegraphics[width=\textwidth]{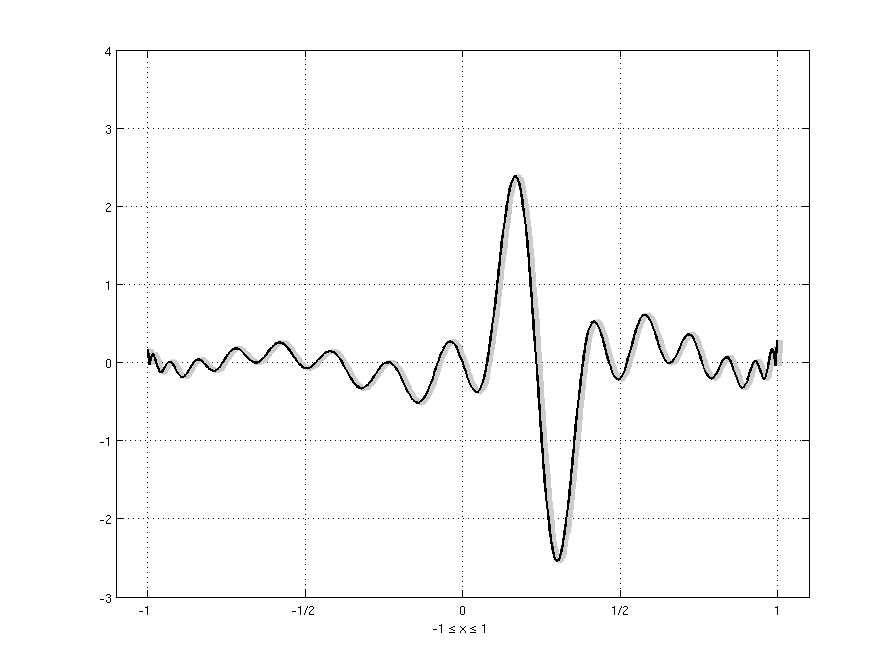}\\
  $\psi_{32,5}^8(x)$. \\
  \includegraphics[width=\textwidth]{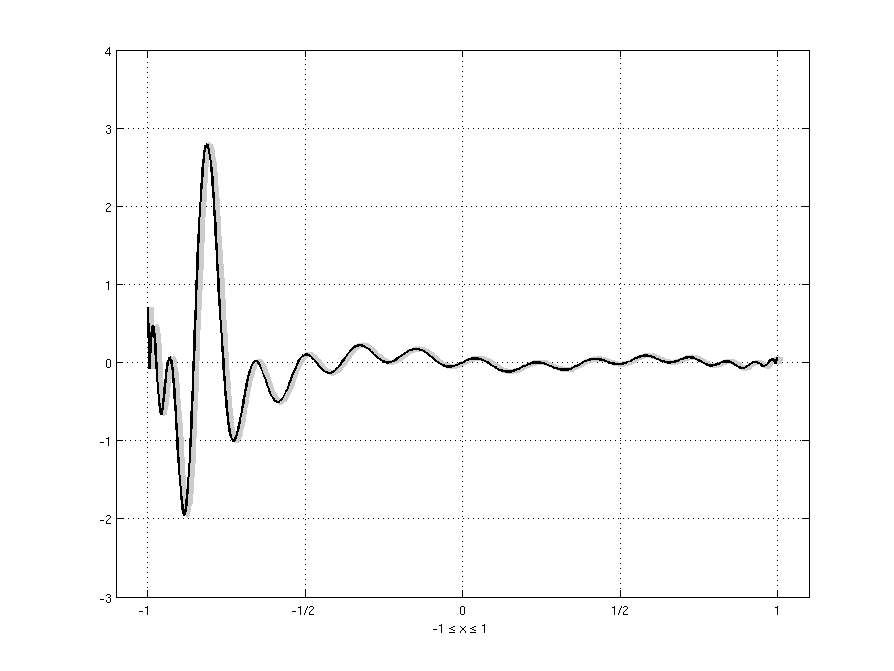}
  \end{minipage}
\caption{Some eigenfunctions $\psi_{n,k}^m$ of the operator $P_n^m M_x P_n^m$ for the Chebyshev polynomials of first kind.}
\label{figure-slepianchebyshev}
\end{figure}

\begin{Exa}
We consider the orthonormal Chebyshev polynomials $t_n$ of first kind defined by
(see \cite[p.~28-29]{Gautschi})
\[ t_0(\cos t) = \ts \frac{1}{\sqrt{\pi}}, \quad t_n(\cos t) = \ts \sqrt{\frac{2}{\pi}} \cos(n t), \quad n \geq 1, \quad \cos t = x.\]
The roots of the Chebyshev polynomials $t_{n+1}$ are given by $x_{n,k} = \cos \frac{2n-2k+3}{2n+2}\pi $, $k = 1, \ldots, n+1$ (see \cite[(6.3.5)]{Szegö}). The normalized
associated polynomials $t_n(x,m)$, $m \geq 1$, correspond to the Chebyshev polynomials $u_n$ of the second kind given by (see \cite[p.~28-29]{Gautschi})
\[ u_n(\cos t) =  \frac{\sin(n+1) t}{\sin t}, \quad n \geq 0.\]
The zeros of the polynomial $u_{n-m+1}$ are given by $x_{n,k}^m = \cos \frac{ n-m+2-k}{n-m+2}\pi $, $k = 1, \ldots, n-m+1$.
Hence, by the formulas \eqref{equation-slepianfunctionsexplicit} and \eqref{equation-slepianfunctionsexplicit0} we get for the eigenfunctions $\psi_{n,k}^m$
the following explicit representation 
\begin{align*}
\psi_{n,k}(\cos t) &= \frac{\kappa_{n,k}}{\pi} \frac{ \cos \frac{n (2n-2k+3)\pi }{2n+2} \cos (n+1)t}{\cos t - \cos \frac{2n-2k+3 }{2n+2}\pi}, \\
\psi_{n,k}^m(\cos t)
&= \frac{\kappa_{n,k}^m}{\sqrt{2\pi}} \frac{ (-1)^{n-m-k+1}\cos (n+1) t +\cos (m-1)t }{\cos t - \cos \frac{n-m+2-k }{n-m+2}\pi}, \quad m \geq 1.
\end{align*}
The constants $\kappa_{n,k}$ can be computed explicitly and are given as (see \cite[Formula (1.1.17)]{Mhaskar})
\[ (\kappa_{n,k})^{-2} = \frac{2n+1+u_{2n}(x_{n,k})}{2\pi}. \]
Some of the eigenfunctions $\psi_{n,k}^m$ are illustrated in Figure \ref{figure-slepianchebyshev}.
\end{Exa}

\section{The localization of the eigenfunctions of $P_n^m M_x P_n^m$}

In this section, we are going to investigate localization properties of the eigenfunctions $\psi_{n,k}^m$.
First of all, we know from \cite[Lemma 2.7]{Erb2012} that the mean value $\eps(P)$ of a polynomial $P(x) = \sum_{l=m}^n c_l p_l(x)$ can be written as
$\eps(P) = \ccc^H \Jj_n^m \ccc$, where $\ccc = (c_m, c_{m+1}, \ldots, c_n)^T$.
A similar characterization can be found for the variance $\vvar(P)$.

\begin{Lem} \label{Lem-characterizationofvariance}
For a normalized polynomial $P(x) = \sum_{l=m}^n c_l p_l(x)$, we have the following characterization of the variance
$\vvar(P)$:
\begin{equation*}
\begin{array}{ll}
\vvar(P) = \ccc^H [\Jj_{n}]^2 \ccc + b_n^2 |c_n|^2 - (\ccc^H \Jj_{n} \ccc)^2,  & \text{if} \quad P \in \pia, \\[2mm]
\vvar(P) = \ccc^H [\Jj_{n}^m]^2 \ccc + b_m^2 |c_m|^2+ b_{n+1}^2 |c_n|^2 - (\ccc^H \Jj_{n}^m \ccc)^2, & \text{if} \quad P \in \pib, \; m \geq 1,
\end{array}
\end{equation*}
with the coefficient vectors $\ccc = (c_m, \ldots, c_n)^T$.
\end{Lem}

\begin{proof}
For $m \geq 1$, we denote by $\mathbf{p}_n^m(x)$ the vector $(p_m(x), \cdots, p_n(x))^H$. Then,
using the three-term recurrence formula \eqref{equation-recursionassociatedsymmetric} and the orthonormality relation of the polynomials $p_l$,
we get for $P(x) = \sum_{l=m}^n c_l p_l(x) \in \pib$, $\|P\|_w = 1$:
\begin{align*}
\vvar(P) &= \int_{-1}^1 \Big|\sum_{l=m}^n c_l x p_l(x)\Big|^2 w(x) dx - \eps(f)^2 \\
&= \int_{-1}^1 \Big|\sum_{l=m}^n c_l \big(b_{l+1} p_{l+1}(x)
+ a_l p_l(x)+ b_{l} p_{l-1}(x) \big)\Big|^2 w(x) dx - \eps(f)^2  \\
&=  \int_{-1}^1 \ccc^H \Jj_n^m {\mathbf{p}}_m^n(x) \cdot {\mathbf{p}}_m^n(x)^H \Jj_n^m \ccc\, w(x) dx +
b_m^2 |c_m|^2 + b_{n+1}^2 |c_n|^2 - \eps(f)^2 \\
&= \ccc^H \Jj_n^m \left( \int_{-1}^1 p_i(x) \overline{p_j(x)} w(x) dx\right)_{i,j = 1}^m \Jj_n^m \ccc + b_m^2 |c_m|^2 + b_{n+1}^2 |c_n|^2 - (\ccc^H \Jj_{n}^m \ccc)^2 \\
&=  \ccc^H [\Jj_{n}^m]^2 \ccc + b_m^2 |c_m|^2 + b_{n+1}^2 |c_n|^2 - (\ccc^H \Jj_{n}^m \ccc)^2.
\end{align*}
For $m=0$, the statement follows analogously but without the term $b_m^2 |c_m|^2$.
\end{proof}

Now, we get the following formulas for the expectation value and the variance of the eigenfunctions
$\psi_{n,k}^m$.

\begin{Lem} \label{Lemma-formulavariances}
For the normalized eigenfunction
$\psi_{n,k}^m$, $1\leq k \leq n-m+1$, corresponding to the eigenvalue $x_{n,k}^m$, we have
\begin{align}
\eps(\psi_{n,k}) = x_{n,k}, \qquad & \vvar(\psi_{n,k}) = b_{n+1}^2 \frac{p_n(x_{n,k})^2}{\sum_{l=0}^n p_l(x_{n,k})^2}, \label{equation-varianceformula1} \\
\eps(\psi_{n,k}^m) = x_{n,k}^m, \qquad & \vvar(\psi_{n,k}^m) = \frac{b_{n+1}^2 p_{n-m}(x_{n,k}^m,m)^2 + b_m^2 }{\sum_{l=0}^{n-m} p_l(x_{n,k}^m,m)^2}. \label{equation-varianceformula2} 
\end{align}
\end{Lem}

\begin{proof}
The statements for the mean value $\eps(\psi_{n,k}^m)$ follow directly from the definition of the $\psi_{n,k}^m$ as eigenfunctions
of the operator $P_n^m M_x P_n^m$. \\
For the variance $\vvar(\psi_{n,k}^m)$ of the normalized eigenfunction $\psi_{n,k}^m$, $m \geq 1$, corresponding to the eigenvalue $x_{n,k}^m$ and with the coefficient vector $\ccc_k$ given in \eqref{equation-slepianfunctionscoefficients}, we can derive from Lemma \ref{Lem-characterizationofvariance} that
\begin{align*}
\vvar(\psi_{n,i}^m) &= \ccc_k^H [\Jj_{n}^m]^2 \ccc_k + b_m^2 |c_{m,k}|^2 + b_{n+1}^2 |c_{n,k}|^2 - (\ccc_k^H \Jj_{n}^m \ccc_k)^2\\
&= (x_{n,k}^m) ^2 (\ccc_k^H \ccc_k)^2 + b_m^2 |c_{m,k}|^2 + b_{n+1}^2 |c_{n,k}|^2 - (x_{n,k}^m \ccc_k^H \ccc_k)^2\\ &= b_m^2 |c_{m,k}|^2 + b_{n+1}^2 |c_{n,k}|^2.
\end{align*}
Inserting the coefficients from \eqref{equation-slepianfunctionscoefficients}, we get the above result. The same argumentation holds also for $m = 0$.
\end{proof}

\begin{Rem}
For the case $m = 0$, the formula \eqref{equation-varianceformula1} for the variance of $\psi_{n,k}$ is a special case of 
a variance formula of the Christoffel-Darboux kernel considered in the proof of \cite[Theorem 2.2]{BreuerLastSimon2010}. 
\end{Rem}

If we want the eigenfunction $\psi_{n,k}^m$ to be localized at the expectation value $x_{n,k}^m$, the variance of $\psi_{n,k}^m$ should be small, especially if $n-m$ gets large.
The question whether the variance in \eqref{equation-varianceformula1} gets small when $n$ is large is linked to a condition known as subexponential growth (see
\cite{BreuerLastSimon2010}, \cite{NevaiZhangTotik1991}). In particular, if the orthonormalization measure $w(x)dx$ is an element of the Nevai class $M(0,1)$, i.e. if the coefficients of the recurrence formula \eqref{equation-recursionorthonormal} attain the limits $\lim_{n\to \infty} a_n = 0$ and $\lim_{n\to \infty} b_n = \frac{1}{2}$, it 
is proven in \cite{NevaiZhangTotik1991} that $\var(\psi_{n,k})$ tends to zero as $n \to \infty$.  
If we restrict the measure $w(x) dx$ to a particular subclass of $M(0,1)$, we can also show in the more general case $m \geq 0$ that the variances in Lemma \ref{Lemma-formulavariances} tend to zero as $n \to \infty$.

\begin{Def} \label{definition-subNevaiclass}
By $M^*(0,1)$, we denote the set of all measures $\mu$ with the following properties:
\begin{enumerate}
 \item $\mu$ is in the Nevai class $M(0,1)$, i.e. $\lim_{n\to \infty} a_n = 0$ and $\lim_{n\to \infty} b_n = \frac{1}{2}$,
 \item $\supp \mu = [-1,1]$,
 \item $\sum_{n = 0}^\infty |a_n| + |b_n - \ts \frac{1}{2}| < \infty$,
\end{enumerate}
where $a_n$ and $b_n$ are the coefficients of the three-term recurrence relation \eqref{equation-recursionassociatedsymmetric} corresponding to the measure $\mu$.
\end{Def}

Examples of weight functions lying in the Nevai subclass $M^*(0,1)$ are, for instance, the Jacobi weight functions (see \cite[p. 79-81]{Nevai}.

For a measure $\mu$ and the corresponding family of orthonormal polynomials $(p_l)_{l \in \Nn}$,
we denote by $\mu_m$ the orthonormalizing measure of the associated polynomials $p_l(x,m)$. In particular, the measure $\mu_m$ is normalized such
that $\mu_m([-1,1]) = 1$. For a measure $\mu$ in the Nevai subclass $M^*(0,1)$, we get the following result:

\begin{Lem} \label{Lemma-absolutelycontinuousassociatedmeasure}
If $\mu \in M^*(0,1)$, then also $\mu_m \in M^*(0,1)$. Moreover, the measures $\mu_m$, $m \geq 1$, are all
absolutely continuous on $[-1,1]$, i.e. $ d\mu_m = w_m dx$.
\end{Lem}

\begin{proof}
Since the coefficients of the three-term recurrence relation \eqref{equation-recursionassociatedsymmetric} of the associated polynomials $p_l(x,m)$ are defined by shifting the
corresponding coefficients of the polynomials $p_l$, the conditions $(1)$ and $(3)$ of Definition \ref{definition-subNevaiclass} are obviously satisfied by the measure $\mu_m$.
The true interval of orthogonality of the sequence of associated polynomials $p_l(x,m)$ is included in the true interval of orthogonality of the original polynomials $p_l(x)$ (see
\cite[Corollary on page 87]{Chihara}). Therefore, $\supp \mu_m \subset \supp \mu = [-1,1]$. Since $\mu_m \in M(0,1)$ is in the Nevai class, also $[-1,1] \subset \supp \mu_m$ holds 
(cf. \cite[Chapter 3.3, Lemma 6]{Nevai})
and, thus, also the property $(2)$ is satisfied. \\
To prove the absolute continuity of $\mu_m$ we use a result of Nevai \cite[Chapter 7, Theorem 40]{Nevai}. This result implies that
if $\mu \in M^*(0,1)$, then the measure $\mu$ consists of an absolutely continuous part $w(x)dx$ on $[-1,1]$ and a point mass $a \delta_{-1} + b \delta_{1}$ on the boundary of $[-1,1]$. Hence, it remains to show that for the associated measures
$\mu_m$, $m\geq 1$ the discrete part vanishes. It is enough if we give the proof for the left hand boundary $x = -1$. In this case, $a=0$ is equivalent to the divergence of the sum $\sum_{l=0}^\infty p_l(-1,m)^2$ (cf. \cite[Theorem 2.1]{Freud}).
By a technique involving chain sequences, Chihara
\cite[Formula (2.18)]{Chihara1989} proved that there is a constant $C_m$ such that
\begin{equation}
 \label{equation-inequalityChihara}
|p_n(-1,m+1)|^2 \geq C_m |p_{n+1}(-1,m)|^2.
\end{equation}
Hence, by a standard induction argument it follows that $\sum_{l=0}^\infty p_l(-1,m)^2$, $m \geq 1 $ diverges, if $\sum_{l=0}^\infty p_l(-1,1)^2$ diverges. So, to complete the proof we have to show the divergence of $\sum_{l=0}^\infty p_l(-1,1)^2$. If $\mu$ is continuous at $x = -1$, then
$\sum_{l=0}^\infty p_l(-1)^2$ diverges, and by \eqref{equation-inequalityChihara} also $\sum_{l=0}^\infty p_l(-1,1)^2$ diverges.
If $\mu$ has a point mass at $x = -1$, then by another result of Chihara \cite[Theorem 3]{Chihara1957}, the measure $\mu_1$ cannot have a point mass at $x = -1$. Hence, in this case the sum $\sum_{l=0}^\infty p_l(-1,1)^2$ also diverges.
\end{proof}

\begin{Thm} \label{Theorem-localizedeigenfunctions}
If the weight function $w$ is in the class $M^*(0,1)$, then
\[\lim_{n\to \infty} \vvar(\psi_{n,k}) = 0, \qquad \lim_{n\to\infty} \vvar(\psi_{n,k}^m) =  0, \qquad m \in \Nn, \]
uniformly for all $k$.
\end{Thm}

\begin{proof}
By Lemma \ref{Lemma-absolutelycontinuousassociatedmeasure}, the measures $w_m(x)dx$ lie in the subclass $M^*(0,1)$, hence also in the Nevai class $M(0,1)$. Therefore, by a result of Nevai, Totik and Zhang \cite[Theorem 2.1]{NevaiZhangTotik1991} we have
\[\lim_{n \to \infty} \sup_{x \in [-1,1]} \frac{|p_n(x,m)|^2}{\sum_{l=0}^n |p_l(x,m)|^2} = 0.\]
Further, by Lemma \ref{Lemma-absolutelycontinuousassociatedmeasure} the associated measures $d\mu_m(x) = w_m(x) dx$, $m \geq 1$, are absolutely continuous on $[-1,1]$. Hence, by \cite[II, Theorem 2.1]{Freud}, also
\[ \lim_{n \to \infty} \frac{1}{\sum_{l=0}^n |p_l(x,m)|^2} = 0\]
uniformly on $[-1,1]$. Therefore, the results of Lemma \ref{Lemma-formulavariances} imply that the variances $\vvar(\psi_{n,k})$ and
$\vvar(\psi_{n,k}^m)$ converge to zero (independently of the choice of $k$) as $n$ tends to infinity.
\end{proof}

\begin{Exa}
For some particular weight functions $w$, it is possible to determine the rate of convergence of the variance $\vvar(\psi_{n,k})$ in Theorem \ref{Theorem-localizedeigenfunctions}. For instance, if the weight $w$ is a generalized Jacobi weight, i.e. if $\supp w = [-1,1]$ and 
\[w(x) = \prod_{i=1}^r (x-t_i)^{\gamma_i}, \qquad -1 = t_1 < t_2 < \cdots < t_{r-1} < t_r = 1, \quad \gamma_i > -1,\]
then the rate of convergence can be determined as (see \cite{Nevai}, Theorem 9.31 and Theorem 6.3.28)
\begin{equation*}
\vvar(\psi_{n,k}) = b_{n+1}^2 \frac{p_n(x_{n,k})^2}{\sum_{l=0}^n p_l(x_{n,k})^2} \sim \frac{\sqrt{1-x_{n,k}^2}}{n}, \qquad 1 \leq k \leq n+1.
\end{equation*}
So, for generalized Jacobi weights, the convergence of $\lim_{n \to \infty} \vvar(\psi_{n,k})$ towards zero is at least linear. The convergence rate is even faster, if 
we choose $k$ such that $x_{n,k}$ is among the $N$ ($N \in \Nn$ fixed) smallest or largest roots of $p_{n+1}(x)$.
\end{Exa}

\section{Approximation of localized functions} \label{section-approximationfunctions}

In this paragraph, we are going to investigate how the decomposition of a bandlimited function $f \in \pib$ in the eigenfunctions $\psi_{n,k}^m$ can be used to approximate functions that are well-localized at a point or a subinterval of $[-1,1]$. In this case, not all of the eigenfunctions $\psi_{n,k}^m$ are needed for a good approximation of the function $f$. We will show that mainly only those eigenfunctions are needed that are located themselves in the region in which $f$ is concentrated. 

From now on we assume that the weight function $w$ lies in the Nevai subclass $M^*(0,1)$. Then, for the Hilbertspace 
\[ L^2([-1,1],w) \ominus \Pi_{m-1} := \overline{\spa \{ p_l:\; l \geq m\}}\]
we can introduce an isometric isomorphism $S_m$ by
\begin{equation}
S_m: \;L^2([-1,1],w) \ominus \Pi_{m-1} \to L^2([-1,1],w_m):\quad S_m p_l(x) := p_{l-m}(x,m), \quad l \geq m.
\end{equation}
If the functions $\phi_{n-m,k}$, $1 \leq k \leq n-m+1$ denote the eigenfunctions of the operator $P_{n-m} M_x P_{n-m}$ on the Hilbert space $L^2([-1,1],w_m)$, we can deduce 
from \eqref{equation-slepianfunctionsfourierseries} that 
\[S_m \psi_{n,k}^m(x) = \phi_{n-m,k}(x) \]
holds. Further, for $\epsilon_m > 0$ we say that a continuous function $f \in L^2([-1,1],w) \ominus \Pi_{m-1}$ is $\epsilon_m$-concentrated on an interval $A \subset [-1,1]$ if
\[ \int_{[-1,1] \setminus A} |S_m f(x)|^2 w_m(x) dx \leq \epsilon_m^2 \|f\|_w^2. \]
An $\epsilon_m$-concentrated function $f$ can be approximated as follows:

\begin{Thm} \label{Theorem-localizedapproximationepsilonconcentrated}
Let $f \in L^2([-1,1],w) \ominus \Pi_{m-1}$ be continuous and $\epsilon_m$-concentrated on the subinterval $A \subset[-1,1]$. Then,
\begin{equation} \label{equation-approximationerror}
\lim_{n \to \infty} \left\| f - \sum_{k:\, x_{n,k}^m \in A} \langle f,\psi_{n,k}^m \rangle_w \psi_{n,k}^m \right\|_w \leq
\epsilon_m \|f\|_w.
\end{equation}
If $A = [\cos \alpha, \cos \beta]$, the number of eigenvalues $x_{n,k}^m$ in $A$ is asymptotically given as
\[ \lim_{n\to\infty} \frac{\# \{k :\; x_{n,k}^m \in A\}}{(n-m)} =  \frac{\alpha-\beta}{\pi}.\]
\end{Thm}

\begin{proof}
We use the isomorphism $S_m$ to shift the error term from the Hilbert space $L^2([-1,1],w) \ominus \Pi_{m-1}$ onto $L^2([-1,1],w_m)$:
\begin{equation} \label{equation-proof41-1} 
\left\| f - \underset{k:\, x_{n,k}^m \in A}{\ds \sum} \langle f,\psi_{n,k}^m \rangle_w \psi_{n,k}^m \right\|_w  
= \left\| S_m f - \sum_{k:\, x_{n,k}^m \in A} \langle S_m f, \phi_{n-m,k}\rangle_{w_m}\phi_{n-m,k} \right\|_{w_m}.
\end{equation}
For an arbitrary $N \in \Nn$, we can assume without restriction that $n$ is large enough such that $N < n-m$. 
By $P_N = \sum_{k=0}^{N} \langle S_m f, p_k(\cdot,m) \rangle_w p_k(\cdot,m)$, we denote the best approximation of $S_m f$ in the subspace $\Pi_N$ of $L^2([-1,1],w_m)$, and
by 
\[E_N(S_m f, w_m) = \inf_{P \in \Pi_N} \|S_m f - P\|_{w_m}  =  \|S_m f - P_N\|_{w_m}\] 
the corresponding error term. Now, using \eqref{equation-proof41-1} and the triangle inequality twice, we get 
\begin{align} \label{equation-estimate1}
& \left| \left\| f - \underset{k:\, x_{n,k}^m \in A}{\ds \sum} \langle f,\psi_{n,k}^m \rangle_w \psi_{n,k}^m \right\|_w  
 - \left\| P_N - \sum_{k:\, x_{n,k}^m \in A} \langle P_N,\phi_{n-m,k} \rangle_{w_m}\phi_{n-m,k} \right\|_{w_m} \right|  \\ 
& \qquad \leq \left\| S_m f - P_N + \sum_{k:\, x_{n,k}^m \in A} \langle P_N - S_m f,\phi_{n-m,k} \rangle_{w_m}\phi_{n-m,k} \right\|_{w_m} \hspace {-3mm} \leq 2 E_N(S_m f, w_m). \notag
\end{align}
From the spectral Theorem \ref{Thm-spectraldecomposition}, we know that the eigenfunctions $\phi_{n-m,k}$ are, up to a normalizing factor,
the fundamental polynomials of Lagrange interpolation with respect to the nodes $x_{n,k}^m$, $1 \leq k \leq n-m+1$. In particular, since $P_N \in \Pi_N \subset \Pi_{n-m}$, we have 
(cf. \cite[Section 3.4]{Szegö})
\[\langle P_N,\phi_{n-m,k} \rangle_{w_m} = \kappa_{n,k}^m P_N (x_{n,k}^m). \]
Hence, if we define the bounded function $g$ on $[-1,1]$ by
\[g(x) := \left\{ \begin{array}{ll} P_N(x) & \text{if $x \in [-1,1] \setminus A$,} \\ 0 & \text{if $x \in A$,} \end{array} \right. \]
then the sum 
\[ \sum_{k:\, x_{n,k}^m \in [-1,1] \setminus A} P_N (x_{n,k}^m) \kappa_{n,k}^m \phi_{n-m,k}\] 
corresponds precisely to the Lagrange interpolant of $g$ at the nodes $x_{n,k}^m$, $1 \leq k \leq n-m+1$.
Therefore, by the Erd\H os-Tur\'an-Theorem (the original result can be found in \cite{ErdoesTuran1937}, in our case we need \cite[Chapter 3, Theorem 2.5]{Freud} with the
parameters $A_n = B_n = 0$) we get in the limit $n \to \infty$:
\begin{align} \label{equation-estimate2}
& \lim_{n \to \infty} \left\| P_N - \sum_{k:\, x_{n,k}^m \in A} \langle P_N,\phi_{n-m,k}\rangle_{w_m} \phi_{n-m,k} \right\|_{w_m}^2 \notag \\ 
& \quad = \lim_{n \to \infty} \left\| \sum_{k:\, x_{n,k}^m \in [-1,1] \setminus A} \langle P_N,\phi_{n-m,k} \rangle_{w_m} \phi_{n-m,k} \right\|_{w_m}^2 \notag \\
& \qquad = \int_{-1}^1 g(x)^2 w_m(x) dx = \int_{[-1,1] \setminus A} P_N(x)^2 w_m(x) dx.
\end{align}
Also by the triangle inequality the following estimate holds:
\begin{align} \label{equation-estimate3}
\left| \left( \int_{[-1,1] \setminus A} P_N(x)^2 w_m(x) dx\right)^{\frac{1}{2}} - \left(\int_{[-1,1] \setminus A} S_m f(x)^2 w_m(x) dx\right)^{\frac{1}{2}} \right|
& \leq E_N(S_m f, w_m).
\end{align}
Combining \eqref{equation-estimate1}, \eqref{equation-estimate2} and \eqref{equation-estimate3}, we can conclude for $n \to \infty$:
\begin{align*} 
& \varlimsup_{n\to \infty} \left| \left\| f - \sum_{k:\, x_{n,k}^m \in A} \langle f,\psi_{n,k}^m \rangle_w \psi_{n,k}^m \right\|_w -  
\left(\int_{[-1,1] \setminus A} S_m f(x)^2 w_m(x) dx\right)^{\frac{1}{2}} \right| \\ 
& \leq \varlimsup_{n\to \infty} \left| \left\| f - \sum_{k:\, x_{n,k}^m \in A} \langle f,\psi_{n,k}^m \rangle_w \psi_{n,k}^m \right\|_w  
- \left\| P_N - \sum_{k:\, x_{n,k}^m \in A} \langle P_N,\phi_{n-m,k} \rangle_{w_m}\phi_{n-m,k} \right\|_{w_m} \right| \\
& \quad + \varlimsup_{n\to \infty} \left| \left\| P_N - \sum_{k:\, x_{n,k}^m \in A} \langle P_N,\phi_{n-m,k} \rangle_{w_m}\phi_{n-m,k} \right\|_{w_m} - 
\left( \int_{[-1,1] \setminus A} P_N(x)^2 w_m(x) dx \right)^{\frac{1}{2}} \right| \\
& \quad + \varlimsup_{n\to \infty} \left| \left( \int_{[-1,1] \setminus A} P_N(x)^2 w_m(x) dx\right)^{\frac{1}{2}} - \left(\int_{[-1,1] \setminus A} S_m f(x)^2 w_m(x) dx\right)^{\frac{1}{2}} \right|
\\ & \leq 3 E_N(S_m f, w_m).
\end{align*}
Since $N$ was choosen arbitrarily, we finally get 
\begin{equation*} 
\lim_{n\to \infty} \left| \left\| f - \sum_{k:\, x_{n,k}^m \in A} \langle f,\psi_{n,k}^m \rangle_w \psi_{n,k}^m \right\|_w -  
\left(\int_{[-1,1] \setminus A} S_m f(x)^2 w_m(x) dx\right)^{\frac{1}{2}} \right| = 0.
\end{equation*}
Inequality \eqref{equation-approximationerror} now follows from the fact that $f$ is $\epsilon_m$-concentrated on $A$. \\
Since the weight function $w$ is in the class $M^*(0,1)$, Lemma \ref{Lemma-absolutelycontinuousassociatedmeasure} ensures that also the associated weight functions $w_m$ are in $M^*(0,1)$. 
This implies $\supp w_m = [-1,1]$ and, by \cite[Theorem 7.29]{Nevai}, that the restricted support of $w_m$ on $[-1,1]$ has measure $2$. 
Therefore, by a well-known result of Erd\H os and Tur\'an (see \cite{ErdoesFreud1974}, \cite{ErdoesTuran1940}) $w_m(x) dx$ is an arc-sine measure which implies
the second statement of Theorem \ref{Theorem-localizedapproximationepsilonconcentrated}.
\end{proof}

\begin{Rem}
The second statement in Theorem 4.1 is not a new result and intended here only as an additional information on the asymptotic number of eigenfunctions involved in 
the approximation process. It is a special case of a general property that for a large class of orthogonal polynomials the asymptotic distribution of the zeros is given by
the arc-sine measure. For weights as the functions $w_m$ this was proven by Erd\H os and Tur\'an in \cite{ErdoesTuran1940}. Far more general conditions leading to the arc-sine property are elaborated in \cite{ErdoesFreud1974}. In particular, it can be shown that every measure in the Nevai class $M(0,1)$ has this property (see \cite[Theorem 5.3]{Nevai}).
\end{Rem}

If a polynomial $P \in \pib$ is localized at the end points $x = -1$ or $x=1$, or if $P$ has a small variance $\vvar(P)$, we obtain the following
error estimates:

\begin{Thm} \label{Theorem-localizedapproximationvariance}
Let $a > 0$ and $I_{-}$ and $I_+$ denote the Intervals $I_{-} = [-1,-1+a]$ and $I_{+} = [1-a,1]$. If $P \in \pib$, $\|P\|_w = 1$, is localized
at the boundary points of $[-1,1]$, we have the following error bounds:
\begin{align}
& \left\| P - \sum_{x_{n,k}^m \in I_{-}} \langle P,\psi_{n,k}^m \rangle_w \psi_{n,k}^m \right\|_w^2 \leq \frac{1+\eps(P)}{a}, \label{equation-errorbound1} \\
& \left\| P - \sum_{x_{n,k}^m \in I_{+}} \langle P,\psi_{n,k}^m \rangle_w \psi_{n,k}^m \right\|_w^2 \leq \frac{1-\eps(P)}{a}. \label{equation-errorbound2}
\end{align}
Further, if $I = [\eps(P)-a, \eps(P) + a] \subseteq [-1,1]$, we get the following error estimate:
\begin{align}
& \left\| P - \sum_{x_{n,k}^m \in I} \langle P,\psi_{n,k}^m \rangle_w \psi_{n,k}^m \right\|_w^2 \leq \frac{\vvar(P)}{a^2}. \label{equation-errorbound3}
\end{align}
\end{Thm}

\begin{proof}
For $P \in \pib$, we have
\begin{align*}
& \left\| P - \sum_{k:\, x_{n,k}^m \in I_{-}} \langle P,\psi_{n,k}^m \rangle_w \psi_{n,k}^m \right\|_w^2 = \sum_{k:\,x_{n,k}^m \in [-1,1] \setminus I_{-}} |\langle P,\psi_{n,k}^m \rangle_w|^2. \\
& \qquad \leq \frac{1}{a} \sum_{k:\,x_{n,k}^m \in [-1,1] \setminus I_{-}} |\langle P,\psi_{n,k}^m \rangle_w|^2 (1 + x_{n,k}^m)
\leq \frac{1}{a} \sum_{k = 1}^{n-m+1} |\langle P,\psi_{n,k}^m \rangle_w|^2 (1 + x_{n,k}^m).
\end{align*}
Since $\|P\|_w^2 = \sum_{k = 1}^{n-m+1} |\langle P,\psi_{n,k}^m \rangle_w|^2 = 1$ and $\sum_{k = 1}^{n-m+1} x_{n,k}^m |\langle P,\psi_{n,k}^m \rangle_w|^2 = \eps(P)$,
we get the stated bound \eqref{equation-errorbound1}. In a similar fashion, the bound \eqref{equation-errorbound2} can be proven. To prove \eqref{equation-errorbound3}, we proceed
also in a simalar way.
\begin{align*}
& \left\| P - \sum_{k:\, x_{n,k}^m \in I} \langle P,\psi_{n,k}^m \rangle_w \psi_{n,k}^m \right\|_w^2
= \sum_{k:\,x_{n,k}^m \in [-1,1] \setminus I} |\langle P,\psi_{n,k}^m \rangle_w|^2 \\
& \quad \leq \frac{1}{a^2} \sum_{k:\,x_{n,k}^m \in [-1,1] \setminus I} |\langle P,\psi_{n,k}^m \rangle_w|^2 (\eps(P) - x_{n,k}^m)^2
\leq \frac{1}{a^2} \sum_{k = 1}^{n-m+1} |\langle P,\psi_{n,k}^m \rangle_w|^2 (\eps(P) - x_{n,k}^m)^2 \\
& \quad = \frac{1}{a^2} \sum_{k = 1}^{n-m+1} |\langle P,\psi_{n,k}^m \rangle_w|^2 ((x_{n,k}^m)^2-\eps(P)^2) \leq \frac{\vvar(P)}{a^2}.
\end{align*}

\end{proof}

\begin{Rem}
Given a normalized polynomial $P \in \pib$, we consider the discrete density function $\rho$ by
\[ \rho(x) = \left\{ \begin{array}{ll} (\langle P, \psi_{n,k}^m \rangle_w)^2 & \text{if $x = x_{n,k}^m$, $k = 1, \ldots n-m+1$}, \\
                                        0 & \text{otherwise}. \end{array} \right. \] 
Then, we can interpret the results of Theorem \ref{Theorem-localizedapproximationvariance} as versions of the Markov and the Chebyshev inequality for a 
$\rho$-distributed random variable. (cf. \cite[p. 114]{Papoulis}).
\end{Rem}

\section{An uncertainty principle for the operators $M_x$ and $P_n^m$}

We are now going to discuss an uncertainty principle related to the operators $M_x$ and $P_n^m$. In particular, we will discuss the trade off between the space localization of
$f$ at the boundary points $x = 1$ and $x = -1$ of $[-1,1]$ and the frequency localization of $f$ in the polynomial subspace $\pib$. The obtained results are very similar to the uncertainty principle stated in 
the theory of Landau, Pollak and Slepian (see \cite{FollandSitaram1997}, \cite{LandauPollak1961}). However, the fact that $M_x$ is not a projection operator will lead to coarser statements and in some extent to differences in the proofs compared to the original setting. A detailed proof of the uncertainty principle in the Landau-Pollak-Slepian theory can be found in \cite[Chapter 2.9]{DymMcKean} and \cite{LandauPollak1961}. An abstract version of the Landau-Pollak-Slepian uncertainty principle involving two arbitrary projection operators on a Hilbert space can be found
in \cite[Part 1, Chapter 3]{HavinJöricke}. An extension of the Landau-Pollak-Slepian uncertainty to more general weight functions is given in \cite{Machluf2008}.

\begin{figure}
\begin{center}
\begin{tikzpicture}[scale=1]

\tikzset{
axis/.style={thick, ->, >=stealth'},
help lines/.style={dashed},
important line/.style={thick},
dot/.style={circle,fill=black,minimum size=4pt,inner sep=0pt,
            outer sep=-1pt}
}

    \def\xmax{4}
    \def\xmin{-3.5}
    
    \fill [fill=orange!20!white] (\xmin, 5) -- (\xmax, 5) -- (5, 0) -- (-5,0) -- cycle;
    \filldraw [fill=orange!10!white] (5, 2.5+0.5*\xmax) arc (0:90: 5 cm - \xmax cm and 0.5 cm ) -- (5,0) -- cycle;
    \filldraw [fill=orange!10!white] (\xmin, 5) arc (90:180: 5 cm + \xmin cm and 0.75 cm ) -- (-5,0) -- cycle;

    \draw[axis] (0,-0.5)  --  (0,6) node(yline)[above] {$\pi_n^m(f)$};
    \draw[axis] (-6,0)  --  (6,0) node(xline)[right] {$\eps(f)$};

    \draw (5,0)  node[dot, label=below: $1$]{};
    \draw (-5,0) node[dot, label=below: $-1$]{};
    \draw (\xmax,0) node[dot, label=below: $x_{n,\max}^m$]{};
    \draw (\xmin,0) node[dot, label=below: $x_{n,\min}^m$]{};
    \draw (5,5) node[dot]{};
    \draw (-5,5) node[dot]{};
    
    \draw (0,5) node[left,xshift=0,yshift=0.3cm]{$1$};

    \draw (-5,0) -- (-5,5) node[anchor=south,xshift=-0.1cm, yshift = 0cm] {$(-1,1)$};
    \draw (5,0) -- (5,5) node[anchor=south,xshift=0.1cm, yshift = 0cm] {$(1,1)$};
    \draw (-5,5) -- (5,5);

    \draw[help lines] (\xmax,0) -- (\xmax,5) node[dot]{};
    \draw[help lines] (\xmin,0) -- (\xmin,5) node[dot]{};
    
    \draw (2,2) node {$A$};
    \draw (4.7,3.5) node {$B_1$};
    \draw (-4.5,3.5) node {$B_2$};
    \draw[<-, >=stealth] (4.8,4.9) -- (5.8,4.5) node[anchor=west] {$C_1$};
    \draw[<-, >=stealth] (-4.7,4.85) -- (-5.8,4.5) node[anchor=east] {$C_2$};

\end{tikzpicture}
\caption{Graphical presentation of the domains $A$, $B_1$, $B_2$, $C_1$ and $C_2$.}
\label{Figure-uncertainty}
\end{center}
\end{figure}
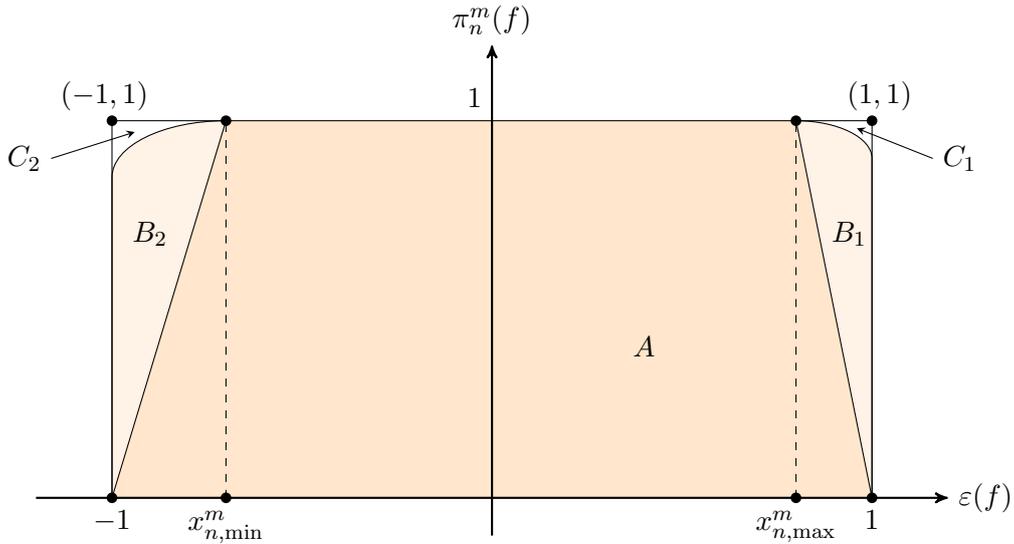

The main results of this section are summarized in Theorem \ref{Theorem-uncertainty} 
and illustrated in Figure \ref{Figure-uncertainty}. The proof of the statements in Theorem \ref{Theorem-uncertainty} is splitted into four lemmas.  
We define
\[ \pi_n^m f := \|P_n^m f\|_w^2 = \sum_{k=m}^n |\langle f, p_k \rangle_w|^2 \]
and start with the first auxiliary result.

\begin{Lem} \label{Lemma-admissibilityoflowerbeta}
Lef $f$, $\|f\|_w = 1$, be a fixed normalized function. Then, for every
$0 \leq \beta \leq \pi_n^m(f)$ there exists a normalized function $g$, $\|g\|_w = 1$, such that $\eps(g) = \eps(f)$ and
$\pi_n^m(g) = \beta$.
\end{Lem}

\begin{proof}
We choose $k > l > n+1$ big enough such
that the three largest eigenvalues $x_1, x_2$ and $x_3$ of the Jacobi matrix $\Jj_k^l$ are larger than $\eps(f)$. This is possible since the weight function $w$ lies in the
class $M^*(0,1)$ and Lemma \ref{Lemma-absolutelycontinuousassociatedmeasure} ensures that also the associated measure $w_l(x) dx \in M^*(0,1)$ is absolutely continuous on $[-1,1]$.
Let $\psi_1, \psi_2$ and $\psi_3$ denote the corresponding eigenfunctions in $\Pi_l^k$. Further, we define $V$ as the $3$-dimensional
vector space spanned by $\psi_1$, $\psi_2$, and $\psi_3$, and $P_V$ as the orthogonal projection operator from $L^2([-1,1],w)$ onto
$V$. Now, we take $\psi$ as a normalized vector in $V$ that is orthogonal to the plane spanned by the vectors $P_V f$ and $P_V M_x f$.
Then, $\eps(\psi) \geq \eps(f)$ and $\langle M_x f, \psi \rangle_w = 0$, $\langle f, \psi \rangle_w = 0$.
In the same way, we construct a normalized vector $\vph \in \Pi_l^k$ with $\eps(\vph) \leq \eps(f)$ and $\langle xf, \vph \rangle_w = \langle f, \vph \rangle_w = 0 $.
Now, since $\eps(f)$ is a continuous functional, by the intermediate value theorem we can find a normalized polynomial $\phi \in \Pi_l^k$ with $\eps(\phi) = \alpha$ and $\langle M_x f, \phi \rangle_w = \langle f, \phi \rangle_w = 0$. Then, we define
\[ g(x) := \sqrt{1-\lambda} f(x) + \sqrt{\lambda} \phi(x), \quad \lambda \in [0,1]. \]
In this way we get a normalized function $g$ with $\|g\|_w = 1$, $\pi_n^m(g) = (1-\lambda) \pi_n^m(f)$ and
\[ \eps(g) = 1-\lambda \eps(f) + \lambda \eps(\phi) = \eps(f). \]
\end{proof}

By $x_{n,\min}^m$ and $x_{n,\max}^m$, we denote the smallest and the largest root of the associated polynomial $p_{n-m+1}(x,m)$. Then, we have as a second
auxiliary result:

\begin{Lem} \label{Lemma-admissibilitybetweenlargestandlowesteigenvalue}
If $x_{n,\min}^m \leq \eps(f) \leq x_{n,\max}^m$, then $\pi_n^m(f)$ can attain all values in the interval $[0,1]$.
\end{Lem}

\begin{proof}
We denote by $\psi_{n,\max}^m$ and $\psi_{n,\min}^m$ the normalized eigenfunctions corresponding to the eigenvalues $x_{n,\max}^m$
and $x_{n,\min}^m$, respectively. Now, for $x_{n,\min}^m \leq \alpha \leq x_{n,\max}^m$, we define the function $f$ by
\[ f = \left(\frac{\alpha-x_{n,\min}^m}{x_{n,\max}^m - x_{n,\min}^m}\right)^{\frac{1}{2}} \psi_{n,\max}^m +
\left(\frac{x_{n,\max}^m-\alpha}{x_{n,\max}^m - x_{n,\min}^m}\right)^{\frac{1}{2}} \psi_{n,\min}^m. \]
Then, $\pi_n^m(f) = \|f\|_w = 1$ and
\[\eps(f) = \frac{\alpha-x_{n,\min}^m}{x_{n,\max}^m - x_{n,\min}^m} x_{n,\max}^m +
 \frac{x_{n,\max}^m-\alpha}{x_{n,\max}^m - x_{n,\min}^m} x_{n,\min}^m = \alpha.\]
Now, Lemma \ref{Lemma-admissibilityoflowerbeta} implies the statement.
\end{proof}

\begin{Lem} \label{Lemma-admissibilitylowerdiagonal}
If  $x_{n,\max}^m \leq \eps(f) < 1$, then $\pi_n^m(f)$ can attain all values in the range $0 \leq \pi_n^m(f) < \frac{ 1 - \eps(f)}{1-x_{n,\max}^m}$.
If  $-1 < \eps(f) \leq x_{n,\min}^m $, then $\pi_n^m(f)$ can attain all values in the range $0 \leq \pi_n^m(f) < \frac{ 1 + \eps(f)}{1+x_{n,\min}^m}$.
\end{Lem}

\begin{proof}
We will prove the statement only for the interval $[x_{n,\max}^m,1)$, the statement for $(-1,x_{n,\max}^m]$ follows by an analagous argumentation.
Since $w(x)dx \in M^*(0,1)$, we can choose as in Lemma \ref{Lemma-admissibilityoflowerbeta} $k>l>n+1$ large enough such that
$1-x_{k,\max}^l < \epsilon$ for an arbitrary $\epsilon >0$. Then, for the eigenfunction $\psi_{k,\max}^l \in \Pi_k^l$ we have $\pi_n^m(\psi_{k,\max}^l) = 0$ and
$1 > \eps(\psi_{k,\max}^l) = x_{k,\max}^l > 1 - \epsilon$.
Now, we define
\[ g(x) = \sqrt{\lambda} \psi_{n,\max}^m(x) +  \sqrt{1-\lambda} \psi_{k,\max}^l(x), \quad \lambda \in [0,1]. \]
Then,
\begin{align*}
 1 - \lambda( 1 - x_{n,\max}^m) > \eps(g) & = \lambda x_{n,\max}^m + (1-\lambda) x_{k,\max}^l > 1-\epsilon - \lambda( 1 - x_{n,\max}^m - \epsilon) \\ & > 1-\epsilon - \lambda( 1 - x_{n,\max}^m),
\end{align*}
and $\pi_n^m(g) = \lambda$. Therefore, we get for $\pi_n^m(g)$:
\[ \frac{1- \eps(g)}{1-x_{n,\max}^m} > \pi_n^m(g) > \frac{1- \eps(g)-\epsilon}{1-x_{n,\max}^m}.\]
Since $\epsilon > 0 $ can be choosen arbitrarily small, we get the desired result from Lemma \ref{Lemma-admissibilityoflowerbeta}.
\end{proof}

Up to now, we showed that most points $(\eps(f), \pi_n^m(f))$ in the rectangle $(-1,1) \times [0,1]$ can be attained for $f\in L^2([-1,1],w)$. However, the next Lemma \ref{Lemma-uncertaintyangle} demonstrates that tuples $(\eps(f), \pi_n^m(f))$ in the upper left and right corner of $(-1,1) \times [0,1]$ are not allowed.

\begin{Lem} \label{Lemma-uncertaintyangle}
If $x_{n,\max}^m \leq \eps(f) < 1$, the values of $\pi_n^m(f)$ are restricted by
\begin{align} \label{equation-uncerataintyangle1a}
\pi_n^m(f)^{\frac{1}{2}} \leq & \frac{(\eps(f)+1)^{\frac{3}{2}} (x_{n,\max}^m+1)^{\frac{1}{2}} + \var(f)^{\frac{1}{2}}(\var(f) + (1+\eps(f))(\eps(f)-x_{n,\max}^m))^{\frac{1}{2}}}{\var(f)+(\eps(f)+1)^2}.
\end{align}
For $ -1 < \eps(f) \leq  x_{n,\min}^m$, the values of $\pi_n^m(f)$ are bounded by
\begin{align} \label{equation-uncerataintyangle2a}
\pi_n^m(f)^{\frac{1}{2}} \leq & \frac{(1-\eps(f))^{\frac{3}{2}} (1-x_{n,\min}^m)^{\frac{1}{2}} + \var(f)^{\frac{1}{2}}(\var(f) + (1-\eps(f))(\eps(f)-x_{n,\min}^m))^{\frac{1}{2}}}{\var(f)+(1-\eps(f))^2}.
\end{align}
A simpler but less accurate upper bound for $\pi_n^m(f)$ is given by
\be \label{equation-uncerataintyangle1b}
\pi_n^m(f) \leq \frac{1}{2} + \frac{1}{2}\Big(\eps(f) x_{n,\max}^m + (1-\eps(f)^2)^{\frac{1}{2}}(1-(x_{n,\max}^m)^2)^{\frac{1}{2}}\Big)
\ee
and
\be \label{equation-uncerataintyangle2b}
\pi_n^m(f) \leq \frac{1}{2} + \frac{1}{2}\Big(\eps(f) x_{n,\min}^m + (1-\eps(f)^2)^{\frac{1}{2}}(1-(x_{n,\min}^m)^2)^{\frac{1}{2}}\Big),
\ee
for $\eps(f)$ in the intervals $[x_{n,\max}^m,1)$ and $(-1,x_{n,\min}^m]$, respectively.
\end{Lem}

\begin{proof}
We will just prove the inequalities \eqref{equation-uncerataintyangle1a} and \eqref{equation-uncerataintyangle1b}. Inequalities \eqref{equation-uncerataintyangle2a} and \eqref{equation-uncerataintyangle2b} follow up to some minor modifications with the same argumentation. Since for $\pi_n^m(f) = 0$
both \eqref{equation-uncerataintyangle1a} and \eqref{equation-uncerataintyangle1b} are satisfied, we will from now on assume that $\pi_n^m(f) > 0$.
Further, we will use the operator $M_{\frac{x+1}{2}}$ on $L^2([-1,1],w)$ defined by $M_{\frac{x+1}{2}} f(x) := \frac{1+x}{2} f(x)$. \\
For a normalized function $f \in L^2([-1,1],w)$ the two functions
$g_1 = \frac{1}{\|M_{\frac{x+1}{2}} f\|_w} M_{\frac{x+1}{2}}f$ and $g_2 = \frac{1}{\|P_n^m f\|_w} P_n^m f$ are also normalized. The sum of the angular distances between the vectors $g_1$ and $f$, and $g_2$ and $f$ is always larger than the angular distance between $g_1$ and $g_2$, i.e.
\be \label{equation-proofuncertainty1} \arccos \Re \langle g_1,f \rangle_w + \arccos \Re \langle g_2,f \rangle_w \geq \arccos \Re \langle g_1,g_2 \rangle_w.\ee
We define the positive selfadjoint operator $M_{\frac{x+1}{2}}^{\frac{1}{2}}$ by $M_{\frac{x+1}{2}}^{\frac{1}{2}} := M_{\sqrt{\frac{x+1}{2}}}$.
Then, for the term $\Re \langle g_1,g_2 \rangle_w$, we can find an upper bound using the Cauchy-Schwarz-inequality and Theorem \ref{Thm-spectraldecomposition}:
\begin{align*}
\Re \langle g_1,g_2 \rangle_w &\leq |\langle g_1,g_2 \rangle_w| = \frac{|\langle M_{\frac{x+1}{2}} f, P_n^m f \rangle_w|}{\|M_{\frac{x+1}{2}} f\|_w \|P_n^m f\|_w} \\ &=
\frac{|\langle M_{\frac{x+1}{2}}^{\frac{1}{2}} f, M_{\frac{x+1}{2}}^{\frac{1}{2}} P_n^m f \rangle_w|}{\|M_{\frac{x+1}{2}} f\|_w \|P_n^m f\|_w}
\leq \frac{ \sqrt{\langle M_{\frac{x+1}{2}} f, f \rangle_w} \sqrt{\langle M_{\frac{x+1}{2}} P_n^m f, P_n^m f \rangle_w} }{\|M_{\frac{x+1}{2}} f\|_w \|P_n^m f\|_w}  \\ &\leq
\frac{ \sqrt{\langle M_{\frac{x+1}{2}} f, f \rangle_w} \sqrt{\frac{x_{n,\max}^m+1}{2}} \sqrt{\langle P_n^m f, P_n^m f \rangle_w} }{\|M_{\frac{x+1}{2}} f\|_w \|P_n^m f\|_w}
= \frac{ \sqrt{\langle M_{\frac{x+1}{2}} f, f \rangle_w} \sqrt{\frac{x_{n,\max}^m+1}{2}} }{\|M_{\frac{x+1}{2}} f\|_w }.
\end{align*}
Now, if we rewrite the expressions $\langle M_{\frac{x+1}{2}} f, f \rangle_w$ and $\|M_{\frac{x+1}{2}} f\|_w$ in terms of $\eps(f)$ and $\var(f)$, we get
\begin{align*}
\Re \langle g_1,g_2 \rangle_w &\leq \frac{ \sqrt{(\eps(f)+1)(x_{n,\max}^m+1)} }{\sqrt{\var(f)+(\eps(f)+1)^2}},\\
\Re \langle g_1,f \rangle_w &= \frac{\eps(f)+1}{\sqrt{\var(f)+(\eps(f)+1)^2}},\\
\Re \langle g_2,f \rangle_w &= \sqrt{\pi_n^m f}.
\end{align*}
Inserting this into inequality \eqref{equation-proofuncertainty1}, we obtain
\be \label{equation-proofuncertainty2} \arccos \frac{\eps(f)+1}{\sqrt{\var(f)+(\eps(f)+1)^2}} + \arccos \sqrt{\pi_n^m f} \geq \arccos \frac{ \sqrt{(\eps(f)+1)(x_{n,\max}^m+1)} }{\sqrt{\var(f)+(\eps(f)+1)^2}}. \ee
Applying the cosine addition formula, this inequality can be rewritten as
\begin{align*}
\sqrt{\pi_n^m f} \leq& \frac{(\eps(f)+1)^{\frac{3}{2}} \sqrt{(x_{n,\max}^m+1)}}{\var(f)+(\eps(f)+1)^2} \\& \qquad + \Big(1-\frac{(\eps(f)+1)^2}{\var(f)+(\eps(f)+1)^2} \Big)^{\frac{1}{2}}
\Big(1-\frac{ (\eps(f)+1)(x_{n,\max}^m+1) }{\var(f)+(\eps(f)+1)^2} \Big)^{\frac{1}{2}} \\
 =& \frac{(\eps(f)+1)^{\frac{3}{2}} (x_{n,\max}^m+1)^{\frac{1}{2}} + \var(f)^{\frac{1}{2}}(\var(f) + (1+\eps(f))(\eps(f)-x_{n,\max}^m))^{\frac{1}{2}}}{\var(f)+(\eps(f)+1)^2}.
\end{align*}
Hence, inequality \eqref{equation-uncerataintyangle1a} is shown. To prove inequality \eqref{equation-uncerataintyangle1b}, we consider inequality
\eqref{equation-proofuncertainty2}. For $0 < a \leq b \leq 1$, the function $\arccos bt - \arccos at$ is a decreasing function of the variable
$t \in [-\frac{1}{b}, \frac{1}{b} ]$. Therefore, if we set
$a = \sqrt{\frac{x_{n,\max}^m+1}{2}} \leq b = \sqrt{\frac{\eps(f)+1}{2}} < 1$ and
\[t = \Big(\frac{2(\eps(f)+1)}{\var(f)+(\eps(f)+1)^2}\Big)^\frac{1}{2} = \Big(\frac{\langle M_{\frac{x+1}{2}} f, f \rangle_w}{\langle M_{\frac{x+1}{2}} f, M_{\frac{x+1}{2}}f \rangle_w}\Big)^\frac{1}{2} \geq 1,\]
we get in inequality \eqref{equation-proofuncertainty2} the upper bound
\[  \arccos \Big(\frac{\eps(f)+1}{2}\Big)^{\frac{1}{2}} + \arccos \sqrt{\pi_n^m f} \geq \arccos \Big(\frac{x_{n,\max}^m+1}{2}\Big)^{\frac{1}{2}}, \]
or equivalently
\be \label{equation-proofuncertainty3}
\sqrt{\pi_n^m f} \leq \frac{1}{2} \Big((\eps(f)+1)^{\frac{1}{2}}(x_{n,\max}^m+1)^{\frac{1}{2}} + (1-\eps(f))^{\frac{1}{2}}(1-x_{n,\max}^m)^{\frac{1}{2}}\Big).
\ee
Taking the square of both sides in \eqref{equation-proofuncertainty3}, we obtain precisely inequality \eqref{equation-uncerataintyangle1b}.
\end{proof}

Now, we introduce the functions $\gamma_1(x)$ and $\gamma_2(x)$ by
\begin{align*}
\gamma_1(x):& [x_{n,\max}^m,1) \to \Rr:\quad \gamma_1(x) := \frac{1}{2} + \frac{1}{2}\Big(x x_{n,\max}^m + (1-x^2)^{\frac{1}{2}}(1-(x_{n,\max}^m)^2)^{\frac{1}{2}}\Big),\\
\gamma_2(x):& (-1,x_{n,\min}^m) \to \Rr:\quad \gamma_2(x) := \frac{1}{2} + \frac{1}{2}\Big(x x_{n,\min}^m + (1-x^2)^{\frac{1}{2}}(1-(x_{n,\min}^m)^2)^{\frac{1}{2}}\Big).
\end{align*}
and the following subdomains of the rectangle $(-1,1) \times [0,1]$ (see Figure \ref{Figure-uncertainty}):

\begin{align*}
A &:= \{(x,y) \in (-1,1) \times [0,1]:\; y < \ts \frac{ 1 - x}{1-x_{n,\max}^m},\;y < \frac{ 1 + x}{1+x_{n,\min}^m}\} \cup \{(x_{n,max}^m,1), (x_{n,\min}^m,1)\}, \\
B_1 &:= \{(x,y) \in (x_{n,max}^m,1) \times [0,1]:\; \ts y \geq \frac{ 1 - x}{1-x_{n,\max}^m},\; y \leq \gamma_1(x)\}, \\
B_2 &:= \{(x,y) \in (-1,x_{n,min}^m) \times [0,1]:\; \ts y \geq \frac{ 1 + x}{1+x_{n,\min}^m},\; y \leq \gamma_2(x)\}, \\
C_1 &:= \{(x,y) \in (x_{n,max}^m,1) \times [0,1]:\; y > \gamma_1(x)\} , \\
C_2 &:= \{(x,y) \in (-1,x_{n,min}^m) \times [0,1]:\; y > \gamma_2(x)\}.
\end{align*}
Finally, we can summarize the results of Lemma \ref{Lemma-admissibilityoflowerbeta}, \ref{Lemma-admissibilitybetweenlargestandlowesteigenvalue}, \ref{Lemma-admissibilitylowerdiagonal}
and \ref{Lemma-uncertaintyangle} as follows.

\begin{Thm} \label{Theorem-uncertainty}
For normalized functions $f \in L^2([-1,1],w)$, all points $(\eps(f),\pi_n^m(f))$ in the domain $A$ can be attained. All points $(\eps(f),\pi_n^m(f))$ in the corners
$C_1$ and $C_2$ cannot be attained.
\end{Thm}

\begin{Rem}
Theorem \ref{Theorem-uncertainty} and its proof based on the Lemmas formulated before are highly inspirated by the uncertainty relation of
the original Landau-Pollak-Slepian theory as described in \cite[Chapter 2.9]{DymMcKean}, \cite[Part 1, Chapter 3]{HavinJöricke}
and \cite{LandauPollak1961}. Lemma \ref{Lemma-admissibilityoflowerbeta} 
reproduces statement F in \cite[Part 1, Section 3.1, p. 95]{HavinJöricke}. However, since $M_x$ is not a projection operator, the proof is altered considerably. 
Lemma \ref{Lemma-admissibilitylowerdiagonal} is an adaption of Case 2 in the proof of \cite[Theorem 2]{LandauPollak1961}. The idea for the proof of Lemma \ref{Lemma-uncertaintyangle}
is taken from \cite[Part 1, Section 3.1 E), p. 95]{HavinJöricke} and the proof of Case 3 in \cite[Theorem 2]{LandauPollak1961}. Due to the fact, that $M_x$ is not a projection operator also here 
the proof differs from the original one. Moreover, the resulting inequalities can not be shown to be sharp. Bounds from below are given in Lemma \ref{Lemma-admissibilitylowerdiagonal}, but it is not yet clear to which extent points $(\eps(f),\pi_n^m(f))$ can be attained in the domains $B_1$ and $B_2$.
\end{Rem}

\end{document}